\documentclass[11pt]{article}
\title{Clique-width of unit interval graphs}
\author{Vadim V. Lozin\thanks{Mathematics Institute, University of Warwick, Coventry CV4 7AL, United Kingdom. 
E-mail: V.Lozin@warwick.ac.uk}}
\date{}
\newtheorem{theorem}{Theorem}

\newtheorem{corollary}[theorem]{Corollary}

\newtheorem{lemma}[theorem]{Lemma}

\newenvironment{proof}[1][Proof]{\textbf{#1.} }{\ \rule{0.5em}{0.5em}}

\begin{document}
\maketitle
\begin{abstract}
The clique-width  is known to be unbounded in the class of unit interval graphs.
In this paper, we show that this is a {\it minimal} hereditary class of unbounded clique-width,
i.e., in {\it every} hereditary subclass of unit interval graphs the clique-width is bounded by a constant.
\end{abstract}

{\it Keywords:} Unit interval graphs; Clique-width

\section{Introduction} 
A graph $G$ is an interval graph if it is the intersection graph of intervals on the real line.
$G$ is a {\it unit} interval graph if all intervals in the intersection model are of the same length.
Unit interval graphs also known in the literature as proper interval graphs \cite{equal} and indifference graphs \cite{Fred}.
These graphs enjoy many attractive properties and find important applications in various fields,
including molecular biology \cite{KSh96}. The structure of unit interval graphs is relatively simple,
allowing efficient algorithms for recognizing and representing these graphs \cite{HShSh01}, as well as for many 
other computational problems \cite{CCC97}.
Nonetheless, some algorithmic problems remain NP-hard when restricted to the class of unit interval graphs
\cite{Marx} and most width parameters are unbounded in this class. 
In the present paper we study the {\it clique-width} of unit interval graphs, which was shown to be unbounded in \cite{GR00}.
Clique-width is a relatively young notion the importance of which is due to the fact that many algorithmic graph
problems which are NP-hard in general become polynomial-time solvable when restricted to graphs of bounded clique-width.
This notion generalizes that of three-width in the sense that graphs of bounded tree-width have bounded clique-width.
The inverse statement is not generally true: there are classes of graphs where the clique-width is bounded but the tree-width 
is not. Cliques (complete graphs) form a trivial example of this type. Notice that every clique is a unit interval graph.
Which other subclasses of unit interval graphs are of bounded clique-width? In the study of this question one may be 
restricted to graph classes which are hereditary in the sense that with any graph they contain all induced subgraphs
of the graph. This restriction is valid due to the fact that the clique-width of a graph cannot be larger than the clique-width of any
of its induced subgraphs \cite{CO00}. Somewhat surprisingly, we show in this paper that the clique-width is bounded
in {\it any} proper hereditary subclass of unit interval graphs. 

We consider simple undirected graphs without loops and multiple edges. For a graph $G$, we denote by $V(G)$ and
$E(G)$ the vertex set and the edge set of $G$ respectively. The neighborhood of a vertex $v\in V(G)$, denote $N_G(v)$, is the 
set of vertices adjacent to $v$. If there is no confusion about $G$ we simply write $N(v)$. We say that $G$ is an $H$-free graph
if no induced subgraph of $G$ is isomorphic to $H$. 
The subgraph of $G$ induced by a subset $U\subseteq V(G)$ will be denoted
$G[U]$. Two vertices of $U$ will be called $U$-similar if they have the same neighborhood outside $U$. 
Clearly, the similarity is an equivalence relation.  The number of equivalence classes of $U$ in $G$ will be denoted 
$\mu_G(U)$ (or simply $\mu(U)$ if no confusion arises). Any subset $U\subseteq V(G)$ with $\mu_G(U)=1$ is called a {\it module} of $G$.
A graph $G$ is said to be {\it prime} if it has no modules $U$ of size $1<|U|<|V(G)|$. When determining the clique-width of graphs in 
a hereditary class $X$ one can be restricted to prime graphs in $X$, because the clique-width of a graph $G$ equals the
clique-width of a maximal prime induced subgraph of $G$ \cite{CO00}.

\section{Canonical unit interval graphs}
\label{sec:can}

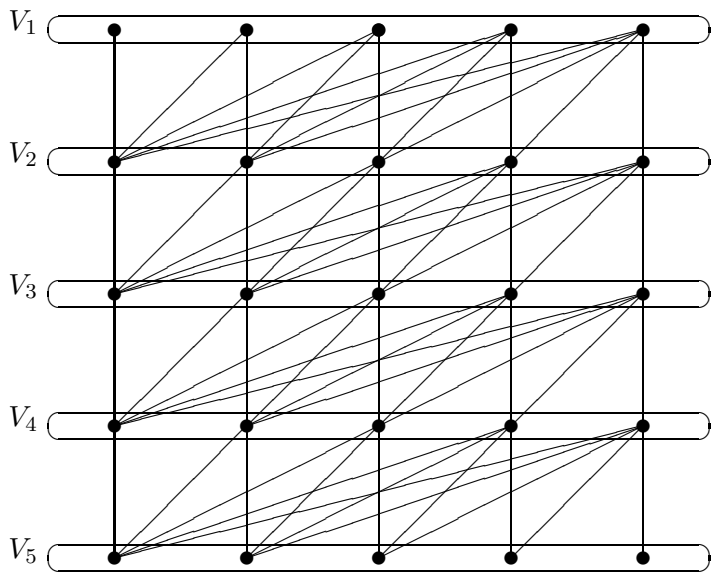
\begin{figure}[ht]
\begin{center}
\begin{picture}(300,220)
\put(50,0){\circle*{5}}
\put(100,0){\circle*{5}}
\put(150,0){\circle*{5}}
\put(200,0){\circle*{5}}
\put(250,0){\circle*{5}}
\put(150,0){\oval(250,10)}
\put(150,50){\oval(250,10)}
\put(150,100){\oval(250,10)}
\put(150,150){\oval(250,10)}
\put(150,200){\oval(250,10)}
\put(10,0){$V_5$}
\put(10,50){$V_4$}
\put(10,100){$V_3$}
\put(10,150){$V_2$}
\put(10,200){$V_1$}
\put(500,50){\circle*{5}}
\put(100,50){\circle*{5}}
\put(150,50){\circle*{5}}
\put(200,50){\circle*{5}}
\put(250,50){\circle*{5}}
\put(50,100){\circle*{5}}
\put(50,50){\circle*{5}}
\put(100,100){\circle*{5}}
\put(150,100){\circle*{5}}
\put(200,100){\circle*{5}}
\put(250,100){\circle*{5}}
\put(50,150){\circle*{5}}
\put(100,150){\circle*{5}}
\put(150,150){\circle*{5}}
\put(200,150){\circle*{5}}
\put(250,150){\circle*{5}}
\put(50,200){\circle*{5}}
\put(100,200){\circle*{5}}
\put(150,200){\circle*{5}}
\put(200,200){\circle*{5}}
\put(150,200){\circle*{5}}
\put(250,200){\circle*{5}}
\put(50,0){\line(0,1){50}}
\put(50,0){\line(1,1){50}}
\put(50,0){\line(2,1){100}}
\put(50,0){\line(3,1){150}}
\put(50,0){\line(4,1){200}}
\put(100,0){\line(0,1){50}}
\put(100,0){\line(1,1){50}}
\put(100,0){\line(2,1){100}}
\put(100,0){\line(3,1){150}}
\put(150,0){\line(0,1){50}}
\put(150,0){\line(1,1){50}}
\put(150,0){\line(2,1){100}}
\put(200,0){\line(0,1){50}}
\put(200,0){\line(1,1){50}}
\put(250,0){\line(0,1){50}}

\put(50,50){\line(0,1){50}}
\put(50,50){\line(1,1){50}}
\put(50,50){\line(2,1){100}}
\put(50,50){\line(3,1){150}}
\put(50,50){\line(4,1){200}}
\put(100,50){\line(0,1){50}}
\put(100,50){\line(1,1){50}}
\put(100,50){\line(2,1){100}}
\put(100,50){\line(3,1){150}}
\put(150,50){\line(0,1){50}}
\put(150,50){\line(1,1){50}}
\put(150,50){\line(2,1){100}}
\put(200,50){\line(0,1){50}}
\put(200,50){\line(1,1){50}}
\put(250,50){\line(0,1){50}}

\put(50,100){\line(0,1){50}}
\put(50,100){\line(1,1){50}}
\put(50,100){\line(2,1){100}}
\put(50,100){\line(3,1){150}}
\put(50,100){\line(4,1){200}}
\put(100,100){\line(0,1){50}}
\put(100,100){\line(1,1){50}}
\put(100,100){\line(2,1){100}}
\put(100,100){\line(3,1){150}}
\put(150,100){\line(0,1){50}}
\put(150,100){\line(1,1){50}}
\put(150,100){\line(2,1){100}}
\put(200,100){\line(0,1){50}}
\put(200,100){\line(1,1){50}}
\put(250,100){\line(0,1){50}}

\put(50,150){\line(0,1){50}}
\put(50,150){\line(1,1){50}}
\put(50,150){\line(2,1){100}}
\put(50,150){\line(3,1){150}}
\put(50,150){\line(4,1){200}}
\put(100,150){\line(0,1){50}}
\put(100,150){\line(1,1){50}}
\put(100,150){\line(2,1){100}}
\put(100,150){\line(3,1){150}}
\put(150,150){\line(0,1){50}}
\put(150,150){\line(1,1){50}}
\put(150,150){\line(2,1){100}}
\put(200,150){\line(0,1){50}}
\put(200,150){\line(1,1){50}}
\put(250,150){\line(0,1){50}}

\end{picture}
\end{center}
\caption{Canonical graph $H_{5,5}$}
\label{fig:H55}
\end{figure}

In this section, we introduce unit interval graphs of a special form that will play an important role
in our considerations. 
Denote by $H_{n,m}$ the graph with $nm$ vertices which can be partitioned into $n$ cliques
$$V_1=\{v_{1,1},\ldots,v_{1,m}\}$$ $$\ldots$$ $$V_n=\{v_{n,1},\ldots,v_{n,m}\}$$ so that for each $i=1,\ldots,n-1$
and for each $j=1,\ldots,m$, vertex $v_{i,j}$ is adjacent to vertices $v_{i+1,1}, v_{i+1,2},\ldots,v_{i+1,j}$
and there are no other edges in the graph. An example of the graph $H_{5,5}$ is given in Figure~\ref{fig:H55}
(for clarity of the picture, each clique $V_i$ is represented by an oval without inside edges).

We will call the vertices of $V_i$ the $i$-th row of $H_{n,m}$,
and the vertices $v_{1,j},\ldots,v_{n,j}$ the $j$-th column of $H_{n,m}$.

It is not difficult to see (and will be clear from the next section) that $H_{n,m}$ is a unit interval graph. Moreover, 
in Section~\ref{sec:univer} we will show that $H_{n,n}$ contains every unit interval graph on $n$ vertices as an induced subgraph.
That's why we will call the graph $H_{n,m}$ {\it canonical} unit interval graph.

Now consider the special case of $H_{n,m}$ when $n=2$. The complement of this graph is bipartite and is known in the literature
under various names such as {\it difference graph} \cite{HPS90} or {\it chain graph} \cite{KKM98}. The latter name is due to the fact that 
the neighborhoods of vertices in each part of the graph form a chain, i.e., the vertices can be ordered under inclusion of their neighborhoods. 
We shall call an ordering $x_1,\ldots,x_k$ {\it increasing} if $i<j$ implies $N(x_{i})\subseteq N(x_j)$ and {\it decreasing} 
if $i<j$ implies $N(x_{j})\subseteq N(x_i)$.
The class of all bipartite chain graphs can be characterized in terms of forbidden induced subgraphs as $2K_2$-free bipartite graphs 
($2K_2$ is the complement of a chordless cycle on $4$ vertices). 
In general, the two parts of a bipartite chain graph can be of different size. But a prime graph in this class 
has equally many vertices in both parts, i.e., it is of the form $H_{2,m}$ with $V_1$ and $V_2$ being independent sets
(see e.g. \cite{FG97}). 

In what follows, we call the complements of bipartite chain graphs {\it co-chain graphs}. Let $G$ be a co-chain graph with
a given bipartition into two cliques $V_1$ and $V_2$, and let $m$ be a maximum number such that $G$ contains the graph $H_{2,m}$ 
as an induced subgraph. Denote by $w_1\in V_1$ and $w_2\in V_2$ two vertices in the same column of $H_{2,m}$ and let 
$W_1:=\{v\in V_1\ | \ N(v)\cap V_2=N(w_1)\cap V_2\}$ 
and $W_2:=\{v\in V_2\ | \ N(v)\cap V_1=N(w_2)\cap V_1\}$. Clearly $W_1\cup W_2$ is a clique and we will call this clique a {\it cluster} of $G$. 
The vertices of $V_1$ that have no neighbors in $V_2$ do not belong to any cluster and we shall call the set of such vertices 
a {\it trivial cluster} of $G$. Similarly, we define a trivial cluster which is a subset of $V_2$. Clearly the set of all
clusters of $G$ defines a partition of $V(G)$.

\section{The structure of unit interval graphs}

To derive a structural characterization of unit interval graphs, we use an ordinary intersection model: 
with each vertex $v$ we associate an interval $I(v)$ on the real line with endpoints $l(v)$ and $r(v)$ such that
$r(v)=l(v)+1$. We will write $I(u)\le I(v)$ to indicate that $l(u)\le l(v)$.    

\begin{theorem}\label{thm:bpg}
A connected graph $G$ is a unit interval graph if and only if the vertex set of $G$ 
can be partitioned into cliques $Q_0,\ldots,Q_t$ in such a way that 
\begin{itemize}
\item[(a)] any two vertices in non-consecutive cliques are non-adjacent,
\item[(b)] any two consecutive cliques $Q_{j-1}$ and $Q_j$ induce a co-chain graph, denoted $G_j$,
\item[(c)] for each $j=1,\ldots,t-1$, there is an ordering of vertices in the clique $Q_j$, 
which is decreasing in $G_{j}$ and increasing in $G_{j+1}$.
\end{itemize}
\end{theorem}

\begin{proof} {\bf Necessity.} Let $G$ be a connected unit interval graph given by an intersection model. 
We denote by $p_0$ a vertex of $G$ with the leftmost interval in the model, i.e., $I(p_0)\le I(v)$ for each vertex $v$.

Define $Q_j$ to be the subset of vertices of distance $j$ from $p_0$ (in the graph-theoretic sense, i.e.,
a shortest path from any vertex of $Q_j$ to $p_0$ consists of $j$ edges). From the intersection model, it is obvious that
if $u$ is not adjacent to $v$ and is closer to $p_0$ in the geometric sense, then it is closer to $p_0$ in the the graph-theoretic sense.
Therefore, each $Q_j$ is a clique.
For each $j>0$, let $p_j$ denote a vertex of $Q_j$ with the rightmost interval in the intersection model. 

We will prove that the partition $Q_0\cup Q_1\cup\ldots\cup Q_t$ satisfies all three conditions of the theorem. 

Condition (a) is due to the definition of the partition. 
Condition (b) will be proved by induction. Moreover, we will show by induction on $j$ that 
\begin{itemize}
\item[(1)] $G_j$ is a co-chain graph,
\item[(2)] $p_{j-1}$ is adjacent to each vertex in $Q_j$,
\item[(3)] for every $v\in Q_i$ with $i\ge j$, $I(p_{j-1})\le I(v)$. 
\end{itemize}
For $j=1$, statements (1), (2), (3) are obvious. To make the inductive step, assume by contradiction that vertices $x_1,x_2\in Q_{j-1}$ and $y_1,y_2\in Q_j$ induce a chordless cycle with edges $x_1y_1$ and $x_2y_2$ (i.e., these vertices induce a $2K_2$ in the complement of $G_j$). 
By the induction hypothesis, both $I(x_1)$ and $I(x_2)$ intersect $I(p_{j-2})$, and also $I(p_{j-2})\le I(y_1), I(y_2)$. Assuming without loss of generality that $I(x_1)\le I(x_2)$,
we must conclude that $I(y_1)$ intersects both $I(x_1)$ and $I(x_2)$, which contradicts the assumption. 
Hence, (1) is correct. 
To prove (2) and (3), consider a vertex $v\in Q_i$, $i\ge j$, non-adjacent to $p_{j-1}$. 
By the induction hypothesis, $I(p_{j-2})$ intersects $I(p_{j-1})$, and also $I(p_{j-2})\le I(v)$, therefore $I(p_{j-1})\le I(v)$, which proves (3). 
Moreover, by the choice of $p_{j-1}$, this also implies that $v$ does not have neighbors in $Q_{j-1}$. 
Therefore, $v\not\in Q_{j}$ and hence (2) is valid.
 
To prove (c), we will show that for every pair of vertices $u$ and $v$ in $Q_j$, 
$N_{G_j}(u)\subset N_{G_j}(v)$ implies $N_{G_{j+1}}(v)\subseteq N_{G_{j+1}}(u)$. Assume the contrary: $s\in N_{G_j}(v) - N_{G_j}(u)$ and $t\in N_{G_{j+1}}(v) - N_{G_{j+1}}(u)$. From (2) we conclude that $s\ne p_{j-1}$. Therefore, $j>1$. Due to the choice of $p_{j-1}$ we have $I(s)\le I(p_{j-1})$, and from (3) we have $I(p_{j-1})\le I(u)$ and $I(p_{j-1})\le I(v)$. Therefore, $I(v)\le I(u)$ by geometric considerations. 
But now, geometric arguments lead us to the conclusion that $tv\in E(G)$ implies $tu\in E(G)$. This contradiction proves (c). 

{\bf Sufficiency.} Consider a graph $G$ with a partition of the vertex set into cliques $Q_0,Q_1,\ldots,Q_t$ satisfying conditions (a), (b), (c). 
We assume that vertices of $$Q_j=\{v_{j,1}, v_{j,2},\ldots, v_{j,k_j}\}$$ are listed in the order that agrees with (c). 
Let us construct an intersection model for $G$ as follows. Each clique $Q_j$ will be represented in the model 
by a set of intervals in such a way that $l(v_{j,i})< l(v_{j,k})< r(v_{j,i})$  whenever $i<k$. 
For $j=0$, there are no other restrictions. For $j>0$, we proceed inductively: for every vertex $u\in Q_j$ with neighbors $v_{j-1,s}, v_{j-1,s+1},\ldots,v_{j-1,k_{j-1}}$ in $Q_{j-1}$, 
we place $l(u)$ between $l(v_{j-1,s})$ and $l(v_{j-1,s+1})$ (or simply to the right of $l(v_{j-1,s})$ if $v_{j-1,s+1}$ does not exist). 
It is not difficult to see that the constructed model represents the graph $G$. 
\end{proof}

\medskip
From this theorem it follows in particular that $H_{n,m}$ is a unit interval graph.  
Any partition of a connected unit interval graph $G$ agreeing with (a), (b) and (c) will be called a {\it canonical partition} of $G$
and the cliques $Q_0,\ldots,Q_t$ the {\it layers} of the partition; cliques $Q_0$ and $Q_t$ will be called marginal layers. 
Any cluster of any co-chain graph $G_j$ in the canonical partition of $G$ will be also called a cluster of $G$. 
\section{$H_{n,n}$ is an $n$-universal unit interval graph}
\label{sec:univer}
The purpose of this section is to show that every unit interval graph with $n$ vertices is contained in the graph $H_{n,n}$ as an induced
subgraph. The proof will be given by induction and we start with the basis of the induction.

\begin{lemma}\label{thm:basis}
The graph $H_{2,n}$ is an $n$-universal co-chain graph.
\end{lemma}

\begin{proof}
Let $G$ be an $n$-vertex co-chain graph with a bipartition into cliques $V_1$ and $V_2$.
We will assume that the vertices of $V_1$ are ordered increasingly according to their neighborhoods in $V_2$,
while the vertices of $V_2$ are ordered decreasingly. The graph $H_{2,n}$ containing
$G$ will be created by adding to $G$
some new vertices and edges. 
Let $W^1,\ldots,W^p$ be the clusters of $G$ and $W_i^j=V_i\cap W^j$. 

For each $W_1^j$ we add to $G$ a set $U_2^j$ of new vertices of size $k=|W_1^j|$ and create on $W_1^j\cup U_2^j$
the graph $H_{2,k}$. Also, create a clique on the set $V'_2=U_2^1\cup W_2^1\cup \ldots\cup U_2^p\cup W_2^p$, and for each $i<j$ 
connect every vertex of $W_1^j$ to every vertex of $U_2^i$.   
Symmetrically, for each $W_2^j$ we add to $G$ a set $U_1^j$ of new vertices of size $k=|W_2^j|$ and create on $W_2^j\cup U_1^j$
the graph $H_{2,k}$. Also, create a clique on the set $V'_1=W_1^1\cup U_1^1\cup\ldots\cup W_1^p\cup U_1^p$, and for each $i<j$ 
connect every vertex of $U_1^j$ to every vertex of $W_2^i\cup U_2^i$. It is not difficult to see that the set $V'_1\cup V'_2$
induces the graph $H_{2,n}$ and this graph contains $G$ as an induced subgraph.
\end{proof}

\medskip
Now we proceed to the general case and assume that every connected unit interval graph $G$ is given together with a canonical partition
$Q_1,\ldots,Q_p$.

\begin{theorem}\label{thm:uni}
Graph $H_{n,n}$ is an $n$-universal unit interval graph.
\end{theorem}

\begin{proof}
Let $G$ be an $n$-vertex unit interval graph. The proof will be given by induction on the number of connected components of $G$.

Assume first that $G$ is connected. We will show by induction on
the number of layers in the canonical partition of $G$ that $H_{n,n}$ contains $G$ as an
induced subgraph, moreover, the $i$-th layer $Q_i$ of $G$ belongs to the
$i$-th row $V_i$ of $H_{n,n}$. The basis of the induction is
established in Lemma~\ref{thm:basis}. Now assume that the
theorem is valid for any connected unit interval graph
with $k\ge 2$ layers, and let $G$ contain $k+1\le n$ layers. For
$j=1,\ldots,k+1$, let $n_j=|Q_j|$ and let $m=n_1+\ldots+n_k$.

Let $H_{k,m}$ be a canonical graph containing the first $k$ layers of $G$ as an induced subgraph.
Now we create an
auxiliary graph $H'$ out of $H_{k,m}$ by
\begin{itemize}
\item[(1)] adding to $H_{k,m}$ the clique $Q_{k+1}$,
\item[(2)] connecting the vertices of $Q_k$ (belonging to $V_k$)
to the vertices of $Q_{k+1}$ as in $G$, 
\item[(3)] connecting the vertices of $V_k-Q_k$ to the vertices of $Q_{k+1}$ so as to make
the existing order of vertices in $V_k$ decreasing in the
subgraph induced by $V_k$ and $Q_{k+1}$. More formally, whenever
vertex $w_{k,i}\in V_k-Q_k$ is connected to a vertex $v\in Q_{k+1}$, every vertex $w_{k,j}$ with $j<i$ must be connected to
$v$ too.
\end{itemize}
According to (2) and (3) the subgraph of $H'$ induced by $V_k$ and
$Q_{k+1}$ is a co-chain graph. We denote this subgraph by $G'$.
Clearly $H'$ contains $G$ as an induced subgraph. To extend $H'$
to a canonical graph containing $G$ we apply the induction
hypothesis twice. First, we extend $G'$ to a canonical co-chain graph
as described in Lemma~\ref{thm:basis}. This will add $m$ new
vertices to the $k+1$-th and $n_{k+1}$ new vertices to $k$-th row of
the graph. Then we use the induction once more to extend the first $k$ rows to a canonical
form. The resulting graph has $k+1\le n$ rows with $n$ vertices
in each row. This completes the proof of the case when $G$ is
connected.

Now assume that $G$ is disconnected. Denote by $G_1$ a connected component of $G$ and by $G_2$ the rest of the graph. 
Also let $k_1=|V(G_1)|$ and $k_2=|V(G_2)|$. Intersection of the first $k_1$ columns and the first $k_1$ rows of $H_{n,n}$ 
induce the graph $H_{k_1,k_1}$, which, according to the above discussion,  contains $G_1$ as an induced subgraph. 
The remaining $k_2$ columns and $k_2$ rows of $H_{n,n}$ induce the graph $H_{k_2,k_2}$, which contains $G_2$ according to the induction hypothesis. 
Notice that no vertex of the $H_{k_1,k_1}$ is adjacent to a vertex of the $H_{k_2,k_2}$. 
Therefore, $H_{n,n}$ contains $G$ as an induced subgraph and the proof is complete.
\end{proof}

\section{Clique-width in subclasses of unit interval graphs}
In this section, we prove that for any proper hereditary subclass $X$ of unit interval graphs, 
the clique-width of graphs in $X$ is bounded by a constant. Let us first recall 
the definition of clique-width.

The {\em clique-width}
of a graph $G$
is the minimum number of labels needed to construct $G$ by means of
the following four operations:
\begin{itemize}
\item[(i)] Creation of a new vertex $v$ with label $i$ (denoted $i(v)$).
\item[(ii)] Disjoint union of two labeled graphs $G$ and $H$
(denoted $G\oplus H$).
\item[(iii)] Joining by an edge each vertex with label $i$ to each vertex with label $j$
($i\not= j$, denoted $\eta_{i,j}$).
\item[(iv)] Renaming label $i$ to $j$
(denoted $\rho_{i\to j}$).
\end{itemize}

Finding the exact value of the the clique-width of a graph is known to be an NP-hard problem \cite{FRRS06}.
In general, this value can be arbitrarily large. Moreover, it is unbounded in many restricted graph families, including 
unit interval graphs \cite{GR00}. On the other hand, in some specific classes of graphs the clique-width is
bounded by a constant. Consider, for instance, a chordless path $P_5$ on five consecutive vertices $a,b,c,d,e$.
By means of the four operations described above this graphs can be constructed as follows:
$$
\eta_{3,2}(3(e)\oplus\rho_{3\to 2}(\rho_{2\to 1}(\eta_{3,2}
(3(d)\oplus\rho_{3\to 2}(\rho_{2\to 1}(\eta_{3,2}(3(c)\oplus\eta_{2,1}(2(b)\oplus 1(a))))))))).
$$
This construction uses only three different labels. Therefore, the clique-width of $P_5$ is at most 3.
Obviously, in a similar way we can construct any chordless path with at most three labels.
This simple example suggests the main idea for the construction of $H_{k,k}$-free unit interval graphs, which
is based on the following lemma (see the introduction for the notation). 

\begin{lemma}\label{lem:cw}
If the vertices of a graph $G$ can be partitioned into subsets $V_1,V_2,\ldots,V_t$ in such a way that
for every $i$ 
\begin{itemize}
\item the clique-width of $G[V_i]$ is at most $k\ge 2$ and 
\item $\mu(V_i)\le l$ and $\mu(V_1\cup\ldots\cup V_i)\le l$,
\end{itemize}
then the clique-width of $G$ is at most $kl$.
\end{lemma}

\begin{proof}
If $G[V_1]$ can be constructed with at most $k$ labels and $\mu(V_1)\le l$, then $G[V_1]$ can be constructed 
with at most $kl$ different labels in such a way that in the process of construction any two vertices in 
different equivalence classes of $V_1$ have different labels, and by the end of the process any two vertices
in the same equivalence class of $V_1$ have the same label. So, the construction of $G[V_1]$ finishes with at most $l$
different labels corresponding to equivalence classes of $V_1$.

Now assume we have constructed the graph $G_i:=G[V_1\cup\ldots\cup V_{i}]$ with the help of $kl$ different labels 
making sure that the construction finishes with a set $A$ of at most $l$ different labels corresponding 
to equivalence classes of $V_1\cup\ldots\cup V_i$. 
Separately, we construct $G[V_{i+1}]$ with the help of $kl$ different labels and complete the construction
with a set $B$ of at most $l$ different labels corresponding to equivalence classes of $V_{i+1}$. 
We choose the labels so that $A$ and $B$ are disjoint. 
Now we use operations $\oplus$ and $\eta$ to build the graph $G_{i+1}:=G[V_1\cup\ldots\cup V_i\cup V_{i+1}]$ out of $G_i$ and $G[V_{i+1}]$. 
Notice that any two vertices in a same equivalence class of $V_1\cup\ldots\cup V_i$ or $V_{i+1}$ 
belong to a same equivalence class of $V_1\cup\ldots\cup V_i\cup V_{i+1}$. Therefore, the construction of $G_{i+1}$ can be 
completed with a set of at most $l$ different labels corresponding to equivalence classes of the graph.
The conclusion now follows by induction.
\end{proof}

\medskip
This lemma implies in particular that 

\begin{corollary}\label{cor:1}
The clique-width of $H_{s,t}$ is at most 3s.
\end{corollary}
\begin{proof}
To build $H_{s,t}$ we partition it into subsets $V_1,V_2,\ldots,V_t$ by including in $V_i$ 
the vertices of the $i$-th column of $H_{s,t}$. Then the clique-width of 
$G[V_i]$ is at most 3. Trivially, $\mu(V_i)=s$. Also, it is not difficult to see that $\mu(V_1\cup\ldots\cup V_i)=s$.
Therefore, the conclusion follows by Lemma~\ref{lem:cw}.
\end{proof}

\medskip 
Now we prove the key lemma of the paper.

\begin{lemma}\label{lem:main}
For every natural $k$, there is a constant $c(k)$ such that the clique-width of any 
$H_{k,k}$-free unit interval graph $G$ is at most $c(k)$.
\end{lemma} 

\begin{proof}
Without loss of generality, we shall assume that $G$ is prime. In particular, $G$ is connected.
To better understand the global structure of $G$, 
let us associate with it another graph which will be denoted $B(G)$. To define $B(G)$ 
we first partition the vertices of $G$ into layers $Q_1,\ldots,Q_t$ as described in Theorem~\ref{thm:bpg}
and then partition each co-chain graph $G_j$ induced by two consecutive cliques $Q_{j-1},Q_j$ into clusters as described in 
Section~\ref{sec:can}. Without loss of generality we may assume that no $G_j$ contains a trivial cluster. Indeed, if such 
a cluster exists, it contains at most one vertex due to primality of $G$. Each $G_j$ contains at most two trivial clusters.
Therefore, by adding at most two vertices to each layer of $G$, we can extend it to a unit interval graph $G'$
such that $G'$ has no trivial clusters, $G'$ contains $G$ as an induced subgraph and $G'$ is $H_{k,k+2}$-free.

With each cluster of $G$ we associate a vertex of the graph $B(G)$ and connect two vertices of $B(G)$ if and only if
the respective clusters have a non-empty intersection. For instance, $B(H_{n,m})$ is a set of $m$ disjoint paths of length
$n-2$ each.
Clearly the vertices of $B(G)$ representing clusters of the same co-chain graph $G_j$ in the partition of $G$ 
form an independent set and we will call this set a level of $B(G)$. In the proof we will use a graphical representation of $B(G)$
obtained by arranging the vertices of the same level on the same horizontal line (different lines for different levels) according to 
the order of the respective clusters in the canonical partition of $G$. 
From this representation it is obvious that $B(G)$ is a plane graph.

Since $G$ is prime, any two clusters of $G$ have at most one vertex in the intersection. 
Therefore, each edge of $B(G)$ corresponds to a vertex of $G$ (this correspondence can be made one-to-one by 
adding to the two marginal levels of $B(G)$ pendant edges representing the vertices of the two marginal layers of $G$). 
 
Now let us consider any $k$ consecutive layers in the canonical partition of $G$ and denote the subgraph of $G$ induced by these 
layers $G^*$. The respective graph $B(G^*)$ will be denoted $B^*$; it has $k-1$ levels denoted $B_1,\ldots,B_{k-1}$.
Since $G$ (and $G^*$) is $H_{k,k}$-free, the two marginal levels of $B^*$ are connected to each other by 
a set ${\cal P}$ of at most $k-1$ disjoint paths. Denote $s=|{\cal P}|$.
 Without loss of generality we may assume that 
the first path in ${\cal P}$ is formed by the leftmost vertices of $B^*$, while the last one by the rightmost vertices of $B^*$.
The $s$ paths of $\cal P$ cut $B^*$ into $s-1$ stripes, i.e., subgraphs induced by two consecutive paths and all the vertices between them.   
 
Since $s$ is the maximum number of disjoint paths connecting $B_1$ to $B_{k-1}$, by Menger's Theorem (see e.g. \cite{Die05}),
these two levels can be separated from each other by a set $S$ of $s\le k-1$ vertices, containing exactly one vertex in each 
of the paths. To visualize this situation, let us draw a curve $\Omega$ that separates $B_1$ from $B_{k-1}$ and crosses $B^*$ at 
precisely $s$ points (the vertices of $S$; no edge is crossed by or belongs to $\Omega$). We claim that without loss of generality
we may assume that this curve traverses each stripe of $B^*$ ``monotonically", meaning that its ``$y$-coordinate'' 
changes within a stripe either non-increasingly or non-decreasingly.
Indeed, assume $\Omega$ has a ``local maximum" within a stripe, and let $v$ be a vertex (below the curve)
that causes this maximum. Obviously, $v$ does not belong to $B_1$ (since otherwise $B_1$ is not separated from $B_{k-1}$), 
and $v$ must have a neighbor at a higher level within the stripe (since there are no trivial blocks in $G$). 
But then the edge connecting $v$ to that neighbor would 
cross $\Omega$, which is impossible according to the definition of $\Omega$. 

The above discussion allows us to conclude that whenever $\Omega$ separates vertices of the same level within a stripe,
the two resulting sets form ``intervals", i.e., their vertices appear in the representation of $B^*$ consecutively.

Now let us translate the above discussion in terms of the graph $G^*$. The partition of the edges of $B^*$ defined by $\Omega$
results in a respective partition of the vertices of $G^*$ into two parts, say $X$ and $Y$. 
Let $Q_i$ be a layer of $G^*$. As we mentioned before, the vertices of $Q_i$ correspond to the edges between two consecutive levels of $B^*$.
We partition these edges and the respective vertices of $Q_i$ into at most $4s-1$ subsets $Q_{i,1},\ldots,Q_{i,4s-1}$ of three types as follows.
The first type consists of $s$ 1-element subsets corresponding to the edges of the $s$ paths of $\cal P$.
For each such an edge $e$, we form at most two subsets of the second type, each consisting of the edges that have a common vertex with $e$ 
and belong to a same stripe. 
The remaining edges form the third group consisting of at most $s-1$ subsets, each representing the edges of the same stripe. 
Observe that the vertices of each $Q_{i,j}$ form an ``interval", i.e., they are consecutive in $Q_i$. The curve $\Omega$ partitions  
each $Q_{i,j}$ into at most two ``subintervals" corresponding to $X$ and $Y$, respectively. 
We claim that no vertex of $Y$ can distinguish 
the vertices of $Q_{i,j}\cap X$. Assume the contrary: a vertex $y\in Y$ is not adjacent to $x_1\in Q_{i,j}\cap X$ but is 
adjacent to $x_2\in Q_{i,j}\cap X$. Then $y\in Q_{i+1}$, $x_2$ and $y$ belong to a same cluster $U$ of $G_{i+1}$, 
while $x_1$ does not belong to $U$. Let $u$ denote the vertex of $B^*$ representing $U$. 
Also, let $e_{x_1},e_{x_2},e_y$ be the edges of $B^*$ corresponding to vertices $x_1,x_2$, and $y$, respectively.
Since $e_{x_2}$ and $e_{y}$ are incident to $u$ but separated by $\Omega$, vertex $u$ belongs to $\Omega$ and hence to the separator $S$.
Therefore, $u$ belongs to a path from $\cal P$. But then $Q_{i,j}$ is of the second type and therefore $e_{x_1}$ must also
be incident to $u$. This contradicts the fact that $x_1$ does not belong $U$. 
This contradiction shows that any two vertices of the same $Q_{i,j}\cap X$ have the same neighborhood in $Y$.
Therefore, $\mu_{G^*}(X)$ is at most the number of different $Q_{i,j}$s, which is at most $k(4s-1)\le 4k^2-5k$.   
Symmetrically, $\mu_{G^*}(Y)\le 4k^2-5k$.

To complete the proof, we partition $G$ into subsets $V_1,\ldots,V_t$ according to the 
following procedure. Set $i:=1$. If the canonical partition of $G$ consists of less than
$k$ layers, then define $V_i:=V(G)$. Otherwise consider the first $k$ layers of $G$ and 
partition the subgraph induced by these layers into sets $X$ and $Y$
as described above. Denote $V_i:=X$ and repeat the procedure with $G:=G-V_i$ and $i:=i+1$. 
By Corollary~\ref{cor:1} each $V_i$ induces a graph of clique-width at most $3k$, and from the above
discussion we know that $\mu(V_i)\le 4k^2-5k$ and $\mu(V_1\cup\ldots\cup V_i)\le 4k^2-5k$. 
Therefore, by Lemma~\ref{lem:cw} the clique-width of $G$ is at most $12k^3-15k$. With the correction on the possible
existence of trivial clusters, we conclude that the clique-width of $G$ is at most $12k^3+72k^2-36k+96$.   
\end{proof}

\begin{theorem}
Let $X$ be a proper hereditary subclass of unit interval graphs. Then the clique-width of graphs in $X$ is bounded by a constant.
\end{theorem}

\begin{proof}
Since $X$ is hereditary, it admits a characterization in terms of forbidden induced subgraphs.
 Since $X$ is a proper subclass of unit interval graphs, it must exclude at least one unit interval graph.
Let $G$ be such a graph with minimum number of vertices. If $|V(G)|=k$, then $G$ is an induce subgraph of $H_{k,k}$
by Theorem~\ref{thm:uni}.
Therefore, $X$ is a subclass of $H_{k,k}$-free unit interval graphs. But then the clique-width of graphs in $X$
is bounded by a constant by Lemma~\ref{lem:main}.
\end{proof}


\end{document}